\documentclass[10pt,twocolumn,letterpaper]{article}

\usepackage{iftex}
\ifPDFTeX
  \usepackage[T1]{fontenc}
  \usepackage[utf8]{inputenc}
  \usepackage{textcomp}
  \usepackage{times}
  \IfFileExists{english.ldf}{\usepackage[english]{babel}}{}
\else
  \usepackage{fontspec}
  \usepackage{polyglossia}
  \setdefaultlanguage{english}
\fi

\usepackage[top=0.75in,bottom=1in,left=0.625in,right=0.625in,columnsep=0.25in]{geometry}
\usepackage{titlesec}
\usepackage{enumitem}
\usepackage{fancyhdr}
\usepackage{microtype}
\usepackage[hidelinks]{hyperref}

\fancypagestyle{plain}{%
  \fancyhf{}
  \fancyhead[R]{\small\thepage}
  
}

\renewcommand{\thesection}{\Roman{section}}
\titleformat{\section}[block]
  {\centering\normalsize\scshape}
  {\thesection.}{0.5em}{}
\titlespacing*{\section}{0pt}{1.4ex plus .5ex}{0.9ex}

\setlist[enumerate]{leftmargin=1.6em,itemsep=2pt,topsep=3pt,label=\arabic*)}
\setlist[itemize]{leftmargin=1.2em,itemsep=2pt,topsep=3pt}

\makeatletter
\renewcommand{\@maketitle}{%
  \begin{center}
    {\fontsize{22}{26}\selectfont\scshape \@title\par}
    \vskip 0.8em
    {\large \@author\par}
  \end{center}
  \vskip 1.2em}
\makeatother

\usepackage{amsmath,amssymb,amsthm}
\usepackage{booktabs,tabularx,array}
\usepackage{tikz}
\usepackage{caption}
\newtheoremstyle{discorsivo}{6pt}{6pt}{\normalfont}{}{\itshape}{.}{.5em}{}
\theoremstyle{discorsivo}
\newtheorem{osservazione}{Remark}

\renewcommand{\thesubsection}{\thesection.\Alph{subsection}}
\titleformat{\subsection}[runin]{\normalfont\normalsize\bfseries}{\thesubsection.}{.5em}{}[.]
\titlespacing*{\subsection}{0pt}{1.2ex plus .3ex}{.5em}
\hypersetup{pdftitle={Mathematical Reality Capture: from the age of discoveries to mathematical topography and surveys},pdfauthor={D. Corradetti},pdflang={en-GB}}
\title{Mathematical Reality Capture:\\from the age of discoveries to\\mathematical topography and surveys}
\author{D. Corradetti}
\date{}
\usepackage{flushend}
\begin{document}
\maketitle\thispagestyle{plain}
\begin{abstract}
A proof answers a precise question, but mathematical understanding also grows through the relationships that its answer makes visible. A point measurement likewise differs from a photograph or a three-dimensional reconstruction of a building. In this work, we argue that, with the help of AI, mathematics should move from a discipline focused on isolated investigations of mathematical objects towards one capable of mapping them systematically. We therefore introduce the idea of \emph{Mathematical Reality Capture}, inspired by digital surveying and the reconstruction of a scene from several viewpoints. An elementary example shows how comparing algebra and geometry can illuminate the same phenomenon, provided that we distinguish diversity of explanations from validity of proof. From this analogy we derive criteria for building verifiable atlases, an experiment with which to assess their usefulness, and a reflection on the role of mathematicians in interpreting them.
\end{abstract}\par\smallskip
{\small\noindent\textit{Keywords.} Artificial intelligence, mathematical understanding, proof, knowledge mapping, Mathematical Reality Capture.\par}

\section{Introduction and motivations: knowing a territory}
\label{sec:intro}

A map lets us choose a route before we have travelled it. Its reader can compare two roads, recognise a valley, or anticipate where a passage will become difficult. Its value arises from the relationships between different observations: an isolated elevation tells us little; elevations, distances, and connections brought together make a territory navigable. In mathematics, too, there is a difference between knowing a route and understanding the landscape in which it lies. A proof can lead us to a conclusion with complete rigour while leaving many questions about the relationship between that result and the rest of the theory.

This experience is familiar to anyone who has encountered two proofs of the same theorem. The second, when it offers a different idea, can change the meaning of the first. A calculation becomes a symmetry; a numerical coincidence turns out to be the trace of a construction; an assumption that seemed technical acquires a geometric reason. The theorem was already proved. What grows is our ability to situate it, generalise it, and use it. This ability has always belonged to mathematical research, from major works of synthesis to conversations at a blackboard.

Artificial intelligence introduces a further possibility. In restricted domains, different systems have been used to suggest relationships that guide researchers' intuition, search for constructions through programs, and produce geometric proofs \cite{davies,trinh,romera}. These results concern specific tasks; they do not establish that mathematical research as a whole can now be automated. They do, however, allow us to formulate a question: if exploring examples and arguments becomes less costly, how can we turn that exploration into organised knowledge?

The question concerns the quality of what is retained. A very large archive may contain repetitions, results stripped of their necessary assumptions, imprecise translations, or conjectures presented as facts. It may also be correct in every part and yet difficult to understand. A task of integration is therefore required: recognising which observations concern the same object, establishing what each adds, making dependencies explicit, and seeking principles through which scattered results can be read together.

In this work, we call \emph{Mathematical Reality Capture}, abbreviated MRC, a proposal for organising this task. The expression borrows the language of digital surveying: acquiring observations from several positions, aligning them, and reconstructing a model that can be queried. In mathematics, the observations become examples, proofs, counterexamples, and alternative representations. The intended product is a \emph{mathematical atlas}: a collection in which connections are as explicit as results, and the status of our knowledge remains recognisable throughout.

\begin{table*}[t]
\centering
\small
\renewcommand{\arraystretch}{1.2}
\begin{tabularx}{\textwidth}{@{}p{2.65cm}XX@{}}
\toprule
Activity & In surveying a territory & In a mathematical atlas \\
\midrule
Exploring & Identifying routes and places & Seeking examples, problems, and conjectures \\
Aligning & Bringing measurements into common frames & Making translations between descriptions explicit \\
Verifying & Checking measurements and control points & Checking proofs, assumptions, and dependencies \\
Reconstructing & Connecting observations in a model & Recognising shared structures and explanations \\
Interpreting & Choosing a map suited to the purpose & Deciding which relationships help us understand \\
\bottomrule
\end{tabularx}
\caption{\textit{This table summarises the functions of the analogy. The correspondences guide a method of working; they do not identify mathematical objects with physical objects.}}
\label{tab:sinossi}
\end{table*}

The proposed contribution is expository and methodological. The search for connections did not begin with AI, and representing mathematics through networks of proofs is already a subject of reflection \cite{barkeshli}. Here we focus on the transition from several descriptions to their integration: what makes a new perspective useful, how to verify its translation, and how to establish whether the whole actually helps us understand. In Section~\ref{sec:parallasse} we begin with an elementary example; in Section~\ref{sec:atlante} we describe the construction of the atlas; in Section~\ref{sec:valutazione} we discuss a possible experimental assessment. The final two sections consider mathematical surveying as a practice and the questions that remain open.

\section{Parallax: seeing the same fact from two positions}
\label{sec:parallasse}

To recognise what a second perspective can add, it helps to choose a fact simple enough for the ideas to remain in the foreground. The sum of the first odd numbers offers precisely this opportunity. No specialist knowledge is needed to follow the calculation; comparing its meanings nevertheless introduces a problem that recurs in more advanced theories.

\subsection{A square within a sum}

Take the odd numbers in succession: one, three, five, seven. Their successive sums are one, four, nine, sixteen. We recognise the squares of the first positive integers. This regularity suggests a statement concerning every positive integer $n$:
\begin{equation}
1+3+5+\cdots+(2n-1)=n^2.
\label{eq:dispari}
\end{equation}
Equation~\eqref{eq:dispari} summarises infinitely many equalities. Checking some of them, even very many, does not yet prove that they all hold. That step requires an argument valid for an arbitrary value of $n$.

One argument consists in writing the same sum twice, once in increasing order and once in decreasing order. We pair the terms occupying the same position in the two rows: one with $2n-1$, three with $2n-3$, and so on. Each pair sums to $2n$, and there are $n$ pairs. The two rows therefore total $2n^2$; a single row is worth $n^2$. The proof is complete and highlights a symmetry: reversing the order pairs the terms into equal sums.

Now imagine small square tiles. One tile forms a square of side one. Adding three along two adjacent sides gives a square of side two; adding another five gives a square of side three. At each step, the new layer has the shape of an L. To pass from side $n$ to side $n+1$, one row requires $n$ tiles and the other $n+1$, counting the corner only once. The layer therefore contains $2n+1$ tiles.

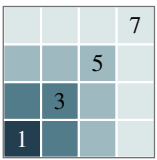
\begin{figure}[htbp]
\centering
\begin{tikzpicture}[x=.5cm,y=.5cm,font=\small]
\definecolor{stratouno}{RGB}{35,62,78}
\definecolor{stratodue}{RGB}{83,124,139}
\definecolor{stratotre}{RGB}{158,186,192}
\definecolor{stratoquattro}{RGB}{220,232,231}
\foreach \x in {0,...,3}{
  \foreach \y in {0,...,3}{
    \pgfmathtruncatemacro{\layer}{max(\x,\y)}
    \ifcase\layer\def\shade{stratouno}\or\def\shade{stratodue}\or\def\shade{stratotre}\or\def\shade{stratoquattro}\fi
    \filldraw[fill=\shade,draw=white,line width=.8pt] (\x,\y) rectangle +(1,1);
  }
}
\draw[black!55,line width=.5pt] (0,0) rectangle (4,4);
\node[text=white] at (.5,.5) {1};
\node at (1.5,1.5) {3};
\node at (2.5,2.5) {5};
\node at (3.5,3.5) {7};
\node[anchor=west,align=left,text width=5.6cm] at (5,2) {Each colour marks one layer.\\[5pt]
The numbers give the tile count\\for the entire layer.\\[5pt]
Four layers form a square\\of sixteen tiles.};
\end{tikzpicture}
\caption{\textit{The sum of the odd numbers viewed as the construction of a square. The figure illustrates four steps; the description of the L-shaped layer extends the argument to every positive integer.}}
\label{fig:quadrato}
\end{figure}

This, too, is a proof. The first square contains one tile, and the count of the layer holds at every step: the construction therefore reaches every positive integer side length. The drawing in Figure~\ref{fig:quadrato} helps us see the mechanism, while the general description ensures that the argument does not depend on the particular square shown. The corresponding algebraic relation is
\begin{equation}
(n+1)^2-n^2=2n+1.
\label{eq:strato}
\end{equation}
In Equation~\eqref{eq:strato}, the difference between consecutive squares counts exactly the tiles in the new layer. Symbol and figure translate into one another.

The first proof suggests pairing the extremes of a progression; the second invites us to study the growth of a figure. On one side, it becomes natural to vary the step of the sequence; on the other, to change the shape or dimension of the construction. The same result thus acquires two families of questions. It is in this opening of possibilities that we recognise a gain in understanding.

\subsection{From physical distance to conceptual difference}

Photogrammetry reconstructs a scene by comparing images taken from different positions. The apparent displacement of a point between images, called \emph{parallax}, helps determine its depth when the geometry of the images is properly controlled \cite{hartley}. Many nearly identical photographs may add less geometric information than a few well-distributed views.

By analogy, we speak of \emph{mathematical parallax} when comparing different descriptions makes visible a relationship that a single description had left in the background. In the odd-number example, the symmetry of pairing and the growth of the square illuminate different aspects of the identity. The analogy concerns our access to the structure rather than the truth of the statement: either proof is already sufficient to establish it.

\begin{osservazione}[Diversity and independence]
Two conceptually different arguments may share axioms, lemmas, and even a substantial part of their reasoning. Their difference does not amount to the probabilistic independence of two measurements. In speaking of parallax, we must specify what changes and what remains shared; the number of proofs does not yield a probability that the theorem is true.
\end{osservazione}

Choosing an elementary example also makes the proposal's limits visible. These two proofs reveal no unknown mechanism and establish no advantage of AI. They offer a readable model of what a system should learn to retain: the steps of the proofs, their point of contact, and the directions in which each suggests proceeding. The interesting problem begins when a theory contains thousands of such encounters, and recognising them requires expertise distributed across different languages.

\section{Building an atlas: acquiring, aligning, explaining}
\label{sec:atlante}

A digital survey becomes usable when its observations find their place within a coherent reconstruction. Individual images may be sharp and still be insufficient: common points must be recognised, distortions checked, and regions actually observed distinguished from those reconstructed under stronger assumptions. In mathematics, the counterpart of this work requires the atlas to preserve the meaning of relationships as well as their presence.

\subsection{Acquiring and retaining the conditions}

The first task is to collect heterogeneous material. For a given result, we may wish to retain examples, counterexamples to its variants, different proofs, equivalent formulations, and unsuccessful attempts. Each item must carry the conditions that make it interpretable: the exact statement, its assumptions, its source, the procedure used, and the status of verification. A regularity observed by computer and a proved theorem occupy different positions in the atlas, even when they express the same formula.

Counterexamples play a particularly important role. If a property ceases to hold when an assumption is removed, the case that refutes it helps draw a boundary. An unsuccessful proof attempt requires greater caution: it may reveal a limitation of the method or of the system searching for it without saying anything about the statement's falsity. Retaining a failure means describing its context, without turning an unsuccessful search into a mathematical impossibility.

This requirement distinguishes the atlas from a simple catalogue of answers. Its reader must be able to trace a result back to the route supporting it. When a calculation depends on a choice of precision or on exploring a finite set, that limit must remain visible. The density of observations acquires value only if we can understand what has actually been observed.

\subsection{Translating without losing structure}

The second task is to put descriptions into correspondence. In the language of surveying, \emph{registration} means aligning observations within common reference frames. In a mathematical atlas, it consists in making translations explicit: which objects correspond, which operations are preserved, and under which assumptions. Saying that a calculation and a drawing tell the same story becomes useful when we can identify what, in the drawing, corresponds to each quantity in the calculation.

For the odd-number identity, the translation is transparent: the integer $n$ becomes a side length measured in tiles, the numerical square becomes a tile count, and the difference becomes the added layer. In a richer theory, correspondences may involve objects very different in appearance. A graph, for example, can be described by a drawing of vertices and connections or by a matrix recording which vertices are adjacent. Moving between these descriptions lets us bring different tools to the same problem.

If a translation is reversible, a natural check consists in moving from one description to the other and then back. Let $T_{ij}$ denote the passage from description $i$ to description $j$. For objects $x$ to which both translations apply, we ask that
\begin{equation}
T_{ji}\bigl(T_{ij}(x)\bigr)\simeq x.
\label{eq:ritorno}
\end{equation}
The sign $\simeq$ means that we recover the same object under the declared equivalence. In the graph example, changing vertex names may change the matrix without changing the connections it represents. This is precisely why the equivalence must be specified before the check.

The return in Equation~\eqref{eq:ritorno} is a consistency check, not a complete guarantee. Two translations could undo the same error and satisfy the return condition while assigning the wrong meaning to operations. Relationships and properties must also be checked. Moreover, some descriptions lose information by construction: the number of vertices alone does not allow a graph to be recovered. In such cases, an inverse cannot be required; what the passage preserves and what it forgets must be declared.

\subsection{From collection to explanation}

The third task is to find a synthesis. The sum of odd numbers offers a miniature example: a long list of equalities is brought under a construction valid for every side length. Someone who understands the L-shaped layer need not remember every numerical result separately. They possess a principle that generates the results and can reconstruct the reason for each one.

We call this search for principles capable of organising many facts \emph{structural compression}. The expression has an affinity with the idea of minimum description length studied in information theory and statistical inference \cite{grunwald}. Here we use it methodologically: brevity alone does not measure mathematical understanding. A sentence may be short because it omits steps, conceals assumptions, or invokes an enormous theorem without explaining it.

A useful synthesis must therefore preserve the ability to reconstruct arguments and address further questions. Its strength becomes apparent when it also clarifies cases not used to formulate it, or reveals the reason for an exception. The atlas should allow its reader to move in both directions: from details to principle and from principle to details. This intellectual reversibility is one of the qualities that distinguish an explanation from a slogan.

\section{What to ask of artificial intelligence}
\label{sec:valutazione}

The surveying analogy becomes a research proposal only if it suggests distinguishable tasks and results that can be checked. A system might produce an atlas rich in intriguing connections without helping anyone solve a new problem. To assess MRC, we must therefore ask what advantage it offers over existing tools and what observations could reveal its limitations.

\subsection{A research cycle guided by gaps}

One possible system would begin with a question and a declared set of sources. It would seek an answer, check its steps, and then try to situate it within its neighbourhood: alternative formulations, weakened assumptions, extreme cases, and arguments based on other tools. The next step would be chosen in response to a recognisable gap. If every proof uses the same lemma, seeking another route might be worthwhile; if a translation works only in simple examples, its points of failure would need to be examined.

Formal libraries provide an important part of this infrastructure. A library such as \emph{mathlib} organises definitions and theorems within a language in which proof steps can be checked by Lean \cite{mathlib}. A generated text and a proof checked by a formal proof assistant nevertheless have different statuses. The check concerns the formalised statement and its assumptions: whether it faithfully expresses the problem we intended to solve still has to be established.

Formalisation can therefore provide solid reference points, while interpretative work connects those points with the atlas's other languages. The comparison with existing programmes becomes concrete: a theorem prover seeks a proof, a program-search system seeks a construction that can be evaluated, and a library makes formalised results reusable. MRC proposes to integrate these functions and assess the integration through the quality of connections and explanations. That advantage remains to be demonstrated.

\subsection{An experiment that could fail}

For an initial test, it is sensible to choose a small, well-documented domain, such as simple undirected graphs with a bounded number of vertices: between two distinct vertices there can be at most one connection, with no assigned direction. Descriptions using connections and matrices, known results, counterexamples, and checked translations would be prepared. Some correspondences would be supplied to the system; others would be reserved for evaluation. Proofs of correctness and exhaustive checks would have to be distinguished from tests performed only on a sample.

The useful comparison would involve one system organised to answer individual questions and another tasked with maintaining the atlas as well. Both would have access to the same sources, the same verification tools, and comparable computational resources. The evaluation would ask how many correct answers they produce, how faithfully they translate objects, and which consequences they can reconstruct in a representation different from the initial one. Test cases would need to be fixed before the experiment.

A second comparison would concern readers. Mathematicians who had not participated in building the atlas could receive explanations from the two systems without knowing how they were produced. Their task would be to address a variant of the problem, identify an indispensable assumption, or recognise an incorrect translation. Times, answers, and reasons would provide more informative evidence than the judgement that an explanation appears elegant. They would still not constitute a universal measure of understanding.

\begin{osservazione}[Testing the proposal]
If, under comparable resources, maintaining the atlas did not improve transfer between representations, error checking, or readers' ability to address new questions, the methodological case for MRC would be weakened. A negative result would be informative: it would indicate where the analogy organises discussion well but does not yet improve practice.
\end{osservazione}

The report should keep coverage, variety of perspectives, quality of translations, verification, and usefulness of syntheses separate. Many nearly identical proofs do not establish great variety; a map rich in objects does not ensure correct translations. A single average could conceal precisely the weakness that the survey should make visible. At this stage, a profile with commentary and examples of success and error is more useful than a ranking of mathematical depth.

\section{Mathematical surveying and the work of interpretation}
\label{sec:rilievo}

A map useful for exploration and one on which a construction can be founded answer different needs. For the first, identifying a promising direction may be enough; for the second, precision, reference points, and conditions of use become decisive. This distinction suggests a further function for the atlas: making explicit what warrants a particular use of the results it contains.

\subsection{Knowing what we rely on}

We call \emph{mathematical surveying} the documented examination of a body of results for a declared purpose. Its subject includes statements, proofs, assumptions, dependencies, and passages between representations. A theorem may be correct and yet fail to apply to the case of interest because a condition is missing or an approximation has changed the problem.

In the odd-number example, the geometric construction counts whole tiles and concerns positive integer side lengths. The algebraic expression can be manipulated in broader contexts, but tile counting does not automatically extend to all of them. A careful survey preserves this distinction: the ease with which we write the same symbols does not warrant transferring their meaning without checking it.

For a complex result, a survey record should make clear which statement has been verified, what assumptions it requires, which dependencies it uses, and which translations have been checked. It should also make excluded cases and open questions visible. Its essential content is argumentative: a sequence of numbers assigned to confidence does not replace an explanation of why a step is justified.

Applications add a further level. Proving a fact about a mathematical model does not by itself establish that the model adequately describes a phenomenon or that the program implementing it is correct. Surveying can clarify the mathematical part of this chain and indicate where other checks begin. The distinction between correctness of deduction, suitability of the model, and reliability of the implementation remains indispensable.

\subsection{Explorers, cartographers, and architects}

If AI tools make some exploration less costly, certain functions may become more prominent. The explorer chooses fruitful questions; the cartographer connects languages and results; the surveyor checks which conditions can be relied upon; the architect seeks principles capable of bringing together what has been discovered. These functions already coexist in a mathematician's work and may alternate within a single afternoon.

Distinguishing them helps us think about the organisation of research without predicting compulsory new professions. A machine may propose thousands of connections, but it remains necessary to decide which are significant, which repeat knowledge already available, and which lead to a problem deserving attention. This judgement includes technical expertise, sensitivity to questions, and the ability to communicate a structure to others.

Teaching, too, might benefit from carefully constructed atlases. In the odd-number identity, one student can enter through calculation and another through the figure; moving between the two makes explicit what both are doing. In more advanced theories, a good translation can serve the same function, offering a point of access without erasing the difficulty. An abundance of routes becomes educational when learners can choose, compare, and justify their own.

A structural limitation remains. Every atlas selects objects, languages, and questions. A well-studied region may appear central because it has more data; a tradition with little formalisation may seem peripheral because its arguments are less accessible to the tools. The map must therefore display its own construction criteria as well. Choosing what is made visible is already part of mathematical interpretation.

\section{Conclusions and future developments}
\label{sec:conclusioni}

In this work, we have developed a reading of AI-assisted mathematics through the surveying analogy: collecting perspectives, clarifying their correspondences, checking their foundations, and seeking a readable synthesis. The small square built from odd numbers shows the most interesting quality of this encounter: a calculation and a figure can give each other meaning, opening different questions from the same truth. MRC proposes to retain such encounters and make them queryable across broader domains. The proposal remains a programme of work: a variety of descriptions does not guarantee an explanation, and automating collection does not guarantee that the connections are correct.

The next reasonable step is to build a bounded atlas, publish its sources and translations, and compare its use with that of an ordinary collection of results. It would be interesting to understand which differences between proofs actually help generalisation, how hidden dependencies can be recognised, and which forms of synthesis remain useful to readers with different backgrounds. A successful atlas should let us return to the proofs, correct a correspondence, and see where our knowledge ends. In this ability to find our bearings, check, and set out again lies the most concrete promise of mathematical surveying: making a wider body of knowledge still inhabitable by thought.

\section*{Use of generative artificial intelligence}

The preparation of this version included the use of a generative artificial intelligence system to rework the source materials for exposition, organise and translate the text, and perform editorial checks. The contribution is conceptual and expository and presents no original experimental results.

\par\medskip
{\footnotesize
\noindent Grupo de F\'\i sica Matem\'atica, Instituto Superior T\'ecnico,
Av. Rovisco Pais, 1049-001 Lisboa, Portugal.\\
\href{mailto:danielecorradetti@tecnico.ulisboa.pt}{\nolinkurl{danielecorradetti@tecnico.ulisboa.pt}}
\par\smallskip
\noindent Departamento de Matematica, Universidade do Algarve,
Campus de Gambelas, 8005-139 Faro, Portugal.\\
\href{mailto:a55499@ualg.pt}{\nolinkurl{a55499@ualg.pt}}\par}
\end{document}